\documentclass{article}
\title{Triangle-degree and triangle-distinct graphs}

\usepackage[total={6in, 9in}]{geometry}
\usepackage{amsthm,amsmath,amssymb}
\usepackage{tikz}
\usepackage{footmisc, soul}

\usepackage{todonotes}
\usepackage{enumerate}

\newtheorem{theorem}{Theorem}
\newtheorem{lemma}[theorem]{Lemma}
\newtheorem{corollary}[theorem]{Corollary}

\newtheorem{pro}[theorem]{Problem}
\theoremstyle{definition}

\renewcommand{\deg}{\operatorname{d}}
\newcommand{\Tri}{\operatorname{t}}

\author{
Zhanar Berikkyzy \thanks{Department of Mathematics, Fairfield University, Fairfield, CT, USA (zberikkyzy@fairfield.edu)} 
\and 
Beth Bjorkman \thanks{Air Force Research Lab, Wright-Patterson Air Force Base, OH, USA (beth.morrison@afrl.af.mil)} 
\and 
Heather Smith Blake \thanks{Department of Mathematics \& Computer Science, Davidson College, Davidson, NC, USA (hsblake@davidson.edu)}
\and 
Sogol Jahanbekam \thanks{Department of Mathematics and Statistics, San Jose State University, San Jose, CA, USA (sogol.jahanbekam@sjsu.edu)}
\and
Lauren Keough \thanks{Department of Mathematics, Grand Valley State University, Allendale, MI, USA (keoulaur@gvsu.edu)} 
\and 
Kevin Moss \thanks{Denver, CO, USA (kmoss646@gmail.com)}
\and 
Danny Rorabaugh \thanks{Alcoa, TN, USA (imnasnainaec@gmail.com)}
\and 
Songling Shan \thanks{Department of Mathematics and Statistics, Auburn University, Auburn, AL, USA (szs0398@auburn.edu)}
}

\begin{document}

\maketitle

\emph{\textbf{Abstract}.}
Let $G$ be a simple graph and $v$ be a vertex of $G$. 
The triangle-degree of $v$ in $G$ is the number of triangles
that contain $v$. While every graph has at least two vertices with the same degree, there are graphs in which every vertex has a distinct triangle-degree. In this paper, 
we construct an infinite family of 
graphs with this property. We also study the vertex degrees and size  
of graphs with this property.

\vspace{2mm}


\vspace{2mm}

\section{Introduction}

In a simple graph $G$, the \emph{triangle-degree} of a vertex $v$, denoted $\Tri_G(v)$, 
is the number of triangles in $G$ that contain $v$.  By the Pigeonhole Principle, every graph has at least two vertices of the same degree, however triangle-distinct graphs exist. These are graphs with at least two vertices in which the triangle-degrees of its vertices are pairwise distinct. Figure~\ref{fig:td7} depicts a triangle-distinct graph with 7 vertices. Erd\H{o}s and Trotter initiated the study of these triangle-distinct graphs by asking for bounds on the number of edges in such a graph \cite{personal},
but very little has been published thus far. Before we describe prior results and our contributions, we first need some notation and terminology. 

In this paper we only consider simple graphs. Let $G$ be a graph with vertex set $V(G)$ and edge set $E(G)$. We denote the cardinality of the vertex set by $n(G)$ and the cardinality of the edge set by $e(G)$. 
For a vertex $v\in V(G)$, the set of neighbors of $v$ 
in $G$ is denoted by $N_G(v)$, and the degree of a vertex $v$ in $G$ is denoted by $d_G(v)=|N_G(v)|$. It is an easy exercise to see that, for any vertex $v$, $\Tri_G(v)$ is the number of edges in $G$ with both endpoints in $N_G(v)$. The closed neighborhood 
of $v$, i.e.,  $N_G(v)\cup \{v\}$, is denoted by $N_G[v]$. 
For $S\subseteq V(G)$, the subgraph of $G$ induced on  $S$ is denoted by $G[S]$, and  $G-S:=G[V(G)\setminus S]$. 
Let $V_1,V_2\subseteq V(G)$ be two disjoint subsets of the vertex set. Then $E_G(V_1,V_2)$ is the set of edges in $G$  with one end in $V_1$ and the other end in $V_2$, and  $e_G(V_1,V_2):=|E_G(V_1,V_2)|$.  We write $E_G(v,V_2)$ and $e_G(v,V_2)$
if $V_1=\{v\}$ is a singleton.   We may drop the subscript $G$, when the graph is understood. The complement of $G$ is denoted $\overline{G}$. 

Nair and Vijayakumar investigated  triangle-degrees in \cite{1}, deriving the value of $\Tri_G(u)+\Tri_{\overline{G}}(u)$ for all graphs. 

\begin{theorem}[Nair and Vijayakumar, \cite{1}]
For any graph $G$ and any $u\in V(G)$,

$$\Tri_G(u)+ \Tri_{\overline{G}}(u)=\left(\sum_{v\in N(u)} d(v)\right)-e(G)+\frac{1}{2}\left(n(G)-d(u)-1\right)\left(n(G)-d(u)-2\right).$$
\end{theorem}

\noindent They also examined the effect of graph composition on triangle-degrees. In particular, given any two graphs $G$ and $H$, the composition graph $G(H)$ has vertex set $V(G)\times V(H)$ where $(u_1,v_1)$ is adjacent to $(u_2,v_2)$ provided $u_1u_2\in E(G)$ or $u_1=u_2$ and $v_1v_2\in E(H)$. 
Nair and Vijayakumar \cite{2} described the triangle-degree for a vertex in the graph $G(H)$ as follows: 

\begin{theorem}[Nair and Vijayakumar, \cite{2}]
Let $G$ and $H$ be two graphs. Then the triangle-degree of any vertex $(u,v)$ in $G(H)$ is given by

\[
\Tri_{G(H)}(u,v)=\Tri_H(v)+e(H)d_G(u)+n(H)d_G(u)d_H(v)+(n(H))^2\Tri_G(u).
\]
\end{theorem}

We are primarily interested in the existence of triangle-distinct graphs and their properties. In the next section, for each $n\geq 7$, we construct a triangle-distinct graph with $n$ vertices.  However, since our construction does not find all triangle-distinct graphs, in Section ~\ref{sec:structural}, 
we establish several properties of triangle-distinct graphs. For example, in a triangle-distinct graph $G$ with $n$ vertices, $\delta(G) \leq n-1-\left(\frac{2n}{3}\right)^{1/3}$ and $\Delta(G)>\sqrt{2n}$ (Theorem~\ref{thm-min-max-degree}). Further, the number of edges in triangle-distinct graphs is between $\frac{\sqrt{2}}{3}n^{3/2}(1-o(1))$ and $\binom{n}{2}-\omega(n)$, where $n$ is the number of vertices (Theorem~\ref{theorem:EdgeLow}).

\section{Construction of Triangle-Distinct Graphs}
\label{sec:Construction}

Similar to the proof that any graph has at least two vertices of the same degree, one can use the Pigeonhole Principle to show there is no triangle-distinct graph on $4$ or fewer vertices. In fact, by a computer-search, the triangle-distinct graph of smallest size is the graph shown in Figure~\ref{fig:td7}.  This graph, $G_7$, has $7$ vertices, $15$ edges, degree sequence $(6, 5, 5, 4, 4, 3, 3)$, and corresponding triangle-degree sequence $(9, 7, 6, 5, 4, 3, 2)$.  The vertices are labeled so that $\Tri(v_1)>\Tri(v_2)> \ldots > \Tri(v_7)$.

\begin{figure}[ht]
	\begin{center}
		\begin{tikzpicture}
		\filldraw (0,0) node[below left]{$v_1$} circle(.08) -- (3,0) node[below right]{$v_6$} circle(.08);
		\filldraw (1,1) node[below]{$v_7$} circle(.08) -- (1,2) node[below right]{$v_5$} circle(.08) -- (1,3) node[left]{$v_4$} circle(.08);
		\filldraw  (2,1) node[below]{$v_3$} circle(.08) -- (2,2) node[above,right]{$v_2$} circle(.08);
		\draw (0,0) to[out=90,in=180] (1,3.5) to[out=0,in=90] (2,2) --(3,0)--(2,1)--(0,0)--(1,1)--(2,1)--(1,3)--(0,0)--(1,2)--(2,2)--(1,3);
		\end{tikzpicture}
		
		\caption{The smallest triangle-distinct graph $G_7$.} \label{fig:td7}
	\end{center}
\end{figure}
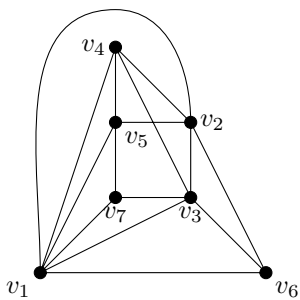

\begin{figure}[ht]
	\begin{center}
		\begin{tikzpicture}
   \draw [red](3,0)--(3,1);
		\filldraw (0,0) node[below left]{$v_1$} circle(.08) -- (3,0) node[below right]{$v_6$} circle(.08);
		\filldraw (1,1) node[below]{$v_7$} circle(.08) -- (1,2) node[below right]{$v_5$} circle(.08) -- (1,3) node[left]{$v_4$} circle(.08);
		\filldraw  (2,1) node[below]{$v_3$} circle(.08) -- (2,2) node[above,right]{$v_2$} circle(.08);
		\draw (0,0) to[out=90,in=180] (1,3.5) to[out=0,in=90] (2,2) --(3,0)--(2,1)--(0,0)--(1,1)--(2,1)--(1,3)--(0,0)--(1,2)--(2,2)--(1,3);
  \filldraw[red] (3,1) node [below right]{$v_8$} circle (.08);
  
            \draw (3,3) node {$G_8$};
            \draw (3,-1.05) node{};
		\end{tikzpicture}
		\hspace{.5in}
  \begin{tikzpicture}
  \draw [red](3,0)--(3,1);
  \foreach \x in {(1,1),(2,2),(2,1),(1,2),(3,1), (3,0),(1,3)} \draw [red](4,1.5)--\x;
            \draw[red] (0,0) to[out=300, in=200] (3,-.75) to [out=20,in=290] (4,1.5);
		\filldraw (0,0) node[below left]{$v_2$} circle(.08) -- (3,0) node[below right]{$v_7$} circle(.08);
		\filldraw (1,1) node[below]{$v_8$} circle(.08) -- (1,2) node[below right]{$v_6$} circle(.08) -- (1,3) node[left]{$v_5$} circle(.08);
		\filldraw  (2,1) node[below]{$v_4$} circle(.08) -- (2,2) node[above,right]{$v_3$} circle(.08);
		\draw (0,0) to[out=90,in=180] (1,3.5) to[out=0,in=90] (2,2) --(3,0)--(2,1)--(0,0)--(1,1)--(2,1)--(1,3)--(0,0)--(1,2)--(2,2)--(1,3);
            \filldraw[red] (3,1) node [below right]{$v_9$} circle (.08);
            \filldraw[red] (4,1.5) node [above right]{$v_1$} circle (.08);
            \draw (3,3) node {$G_9$};
		\end{tikzpicture}
		\caption{The construction of triangle-distinct graphs $G_8$ and $G_9$. The red vertices and edges are the ones that have been added to $G_7$ to create $G_8$ and $G_9$ and vertex labels have been updated as prescribed by the algorithm.} \label{fig:td7}
	\end{center}
\end{figure}
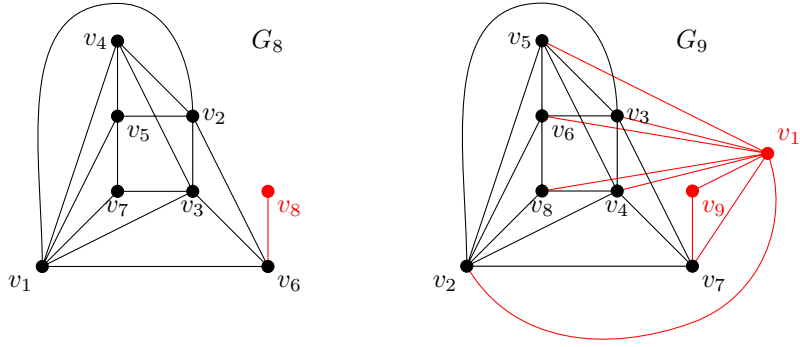

Starting with $G_7$, we give a recursive construction to define a triangle-distinct graph $G_n$ with $n$ vertices for each $n\geq 7$. Note: To simplify notation, for any $v\in V(G_n)$, we write $\deg_n(v)$ and $\Tri_n(v)$ for the degree and triangle-degree, respectively, of $v$ in $G_n$.

The graph $G_7$ is given in Figure~\ref{fig:td7}. Suppose $G_n$ has been defined for some odd $n\geq 7$ with vertices labeled such that $\Tri_n(v_1) > \Tri_n(v_2) > \ldots > \Tri_n(v_n)$. Let $1\leq k\leq n$ be the smallest integer such that $\deg_{n}(v_k)$ equals the minimum degree of $G_n$. To construct $G_{n+1}$, start with a copy of $G_n$, add a vertex labeled $v_{n+1}$ and the edge $v_kv_{n+1}$.  To create $G_{n+2}$,  start with a copy of $G_{n+1}$, relabel each vertex by increasing each subscript by 1, add a single vertex labeled $v_1$ and edges $\{v_1v_i: 2\leq i \leq n+1\}$. 

Now we will use induction to prove that the constructed graphs are in fact triangle-distinct. 

\begin{theorem}
For each $n\geq 7$, $G_n$ is triangle-distinct with vertices labeled such that 
\[\Tri_n(v_1)>\Tri_n(v_2) > \ldots > \Tri_n(v_n).\]
Further, $\deg_n(v_1) \geq \deg_n(v_2) \geq \ldots \geq \deg_n(v_n).$
 If $n$ is odd, we have the additional properties that $\Tri_n(v_n)>0$, $\deg_n(v_n)>0$, $e(G_n) = \Tri_n(v_1) + \deg_n(v_1)$, and $G_n$ is not regular. 
\end{theorem}

\begin{proof}
First consider $G_7$, given in Figure~\ref{fig:td7}, where each of the desired properties can be quickly verified. 

For induction, suppose that for some odd $n\geq 7$, $G_n$ has all of the properties stated in the theorem. We will prove that $G_{n+1}$ and $G_{n+2}$ both satisfy the conditions in the theorem. 

Based on the construction of $G_{n+1}$, we have $\Tri_{n+1}(v_{n+1}) = 0$, $\deg_{n+1}(v_{n+1}) = 1$, $\Tri_{n+1}(v_k) = \Tri_n(v_k)$, and $\deg_{n+1}(v_k) = \deg_n(v_k) + 1$. For each $1\leq i \leq n$ with $i\neq k$, we have $\Tri_{n+1}(v_i) = \Tri_{n}(v_i)$ and $\deg_{n+1}(v_i) = \deg_n(v_i)$. In particular, since $G_n$ is not regular, then $k\neq 1$, so $\Tri_{n+1}(v_1) = \Tri_{n}(v_1)$ and $\deg_{n+1}(v_1) = \deg_n(v_1)$ which will be used in the proof below that $G_{n+2}$ has the desired properties. 

Because $\Tri_n(v_n)>0$, we have 
\[\Tri_{n+1}(v_1) > \Tri_{n+1}(v_{2}) > \ldots > \Tri_{n+1}(v_{n+1})=0.\]
By the choice of $k$ and because $\deg_n(v_n)>0$, we have \[\deg_{n+1}(v_1) \geq \deg_{n+1}(v_2) \geq \ldots \geq \deg_{n+1}(v_{n+1})=1.\]  

Now let's consider $G_{n+2}$. Observe that for each $2\leq i \leq n+2$ we have $\Tri_{n+2}(v_i)= \Tri_{n+1}(v_{i-1}) + \deg_{n+1}(v_{i-1})$. 
Therefore 
\[\Tri_{n+2}(v_2) > \Tri_{n+2}(v_3) > \ldots > \Tri_{n+2}(v_{n+2})=1.\] 
Further, $\Tri_{n+2}(v_1) = e(G_{n+1}) = e(G_n) + 1$ while $\Tri_{n+2}(v_2) = \Tri_{n+1}(v_1) + \deg_{n+1}(v_1)= \Tri_n(v_1) + \deg_n(v_1) = e(G_n)$. 
Therefore, $\Tri_{n+2}(v_1) > \Tri_{n+2}(v_2)$.
 On the other hand, $\deg_{n+2}(v_1) = n+1$, the maximum possible degree, while $\deg_{n+2}(v_i) = \deg_{n+1}(v_i) +1$ for all $2\leq i \leq n+2$. So 
 \[\deg_{n+2}(v_1) \geq \deg_{n+2}(v_2) \geq \ldots \geq \deg_{n+2}(v_{n+2})=2.\]
 Also note that $\Tri_{n+2}(v_1) + \deg_{n+2}(v_1) = e(G_{n+1}) + n+1 = e(G_{n+2})$ as desired. 
 Now $\Tri_{n+2}(v_1) > \Tri_{n+2}(v_{n+2}) = 1$ means that $\Tri_{n+2}(v_1) \geq 2$.
 A vertex with triangle-degree at least 2 has degree at least 3, so $\deg_{n+2}(v_1) \geq 3 > \deg_{n+2}(v_{n+2})$, which means $G_{n+2}$ is not regular.
 This completes the proof.
\end{proof}

\section{Structural properties of triangle-distinct graphs}\label{sec:structural}

In this section, we study the properties of triangle-distinct graphs in terms of 
the minimum and maximum degrees and number of edges.  
For any graph $G$ with $S\subseteq V(G)$, let $\overline{S}$ be $V(G)\setminus S$. 
The following lemma, which relates the degrees and triangle-degrees of a given vertex in a graph and its complement, will be useful in our main results later in this section. 

\begin{lemma} \label{comp}
	Let $G$   be a graph on $n$ vertices and $u \in V(G)$. 
	If $G$ has two distinct vertices $u$ and $v$ with the same degree and \[ \Tri_{\overline{G}}(u) + e_{\overline{G}}\left(N_G(u),\overline{N_{G}[u]}\right) =  \Tri_{\overline{G}}(v) + e_{\overline{G}}\left(N_G(v),\overline{N_{G}[v]}\right),\] then $G$ is not triangle-distinct.
	
\end{lemma}

\begin{proof} Note that $\Tri_{\overline{G}}(u)$ counts the number of edges $xy$
	in $\overline{G}$ with $x,y\not\in N_G(u)$. Counting the pairs of vertices with neither endpoint in $N_G(u)$, we see 
	\[\binom{n-1-d_G(u)}{2} = \left| \left\{ xy\notin E(G): x,y\notin N_G[u] \right\} \right| +  \left| \left\{ xy\in E(G): x,y\notin N_G[u] \right\} \right| = \Tri_{\overline{G}}(u)+e(G-N_G[u]).\]
	Therefore,  
	\begin{align*}
	e(G) &=
	e_G(u,N_G(u))
	+ e(G[N_G(u)])
	+ e_G(N_G(u), \overline{N_G[u]})
	+ e(G-N_G[u]) \\
	&= d_G(u) + \Tri_G(u) + e_G(N_G(u), \overline{N_G[u]}) + \left(\binom{n-1-d_G(u)}{2} - \Tri_{\overline{G}}(u) \right).
	\end{align*}
	Solving the above equation for $\Tri_{G}(u)$, we find,
	\begin{align*}
	 \Tri_G(u)  &= e(G) - d_G(u) - e_G(N_G(u), \overline{N_G[u]}) - \binom{n-1-d_G(u)}{2} + \Tri_{\overline{G}}(u).
	\end{align*}
	Based on the construction of a graph complement, $n-1-d_G(u) = d_{\overline{G}}(u)$. Also, thinking about a complete bipartite graph with one set of size $d_G(u)$ and the other of size $d_{\overline{G}}(u)$, we find \[d_{\overline{G}}(u)d_{G}(u) = e_G(N_G(u), \overline{N_G[u]}) + e_{\overline{G}}\left(N_G(u), \overline{N_G[u]}\right).\]
	Through substitutions, we find: 
	\begin{eqnarray*}
	\Tri_G(u)&=& e(G) - d_{G}(u)-\left(d_{\overline{G}}(u)d_{{G}}(u) - e_{\overline{G}}\left(N_G(u), \overline{N_G[u]}\right)\right)- \binom{d_{\overline{G}}(u)}{2} + \Tri_{\overline{G}}(u)\\
	&=& e(G) - \left(1 + d_{\overline{G}}(u)\right)d_{{G}}(u) - \binom{d_{\overline{G}}(u)}{2} + \Tri_{\overline{G}}(u) + e_{\overline{G}}\left(N_G(u), \overline{N_G[u]}\right)
	\end{eqnarray*}

Now consider distinct vertices $u$ and $v$ which have the same degree and suppose \[\Tri_{\overline{G}}(u) + e_{\overline{G}}\left(N_G(u),\overline{N_{G}[u]}\right) =  \Tri_{\overline{G}}(v) + e_{\overline{G}}\left(N_G(v),\overline{N_{G}[v]}\right)\] as in the hypothesis of the lemma. Since $d_G(u) = d_G(v)$ and consequently $d_{\overline{G}}(u) = d_{\overline{G}}(v)$, then
\[ e(G) - \left(1 + d_{\overline{G}}(u)\right)d_{{G}}(u) - \binom{d_{\overline{G}}(u)}{2} = e(G) - \left(1 + d_{\overline{G}}(v)\right)d_{{G}}(v) - \binom{d_{\overline{G}}(v)}{2}. \]
Thus $\Tri_G(u) = \Tri_G(v)$ so $G$ is not triangle-distinct.	
\end{proof}

\begin{theorem}\label{thm-min-max-degree} 
	Let $G$ be a graph of order $n\geq 3$. Then we have the following:
	\begin{enumerate}[(a)]
		\item If $\Delta(G) \leq \sqrt{2n}$,  
		then $G$ has distinct vertices $u,v \in V(G)$ with $\Tri_G(u)=\Tri_G(v)$.
		\item  If $\delta(G) > n-1-\left(2n/3\right)^{1/3}$, then $G$ has two distinct vertices $u,v \in V(G)$ with $d(u) = d(v)$ and $\Tri_G(u)=\Tri_G(v)$. 
		\item If $G$ is $d$-regular and triangle-distinct, then  \[
		\sqrt{2n} < d \leq n - \sqrt{2n/3}.
		\]
	\end{enumerate}
\end{theorem}

\proof 

For (a), note that for each vertex $v\in V(G)$, there are at most $\binom{\sqrt{2n}}{2}$ edges in $G[N_G(v)]$. Since $\binom{\sqrt{2n}}{2}<n-1$, the triangle-degree for each vertex is between $0$ and $n-2$. 
However, the graph $G$ has $n$ vertices, the Pigeonhole Principle guarantees us two vertices with the same triangle-degree.

For (b), since $\delta(G) > n-1-(2n/3)^{1/3}$, we have $\Delta(\overline{G}) < (2n/3)^{1/3}$. Because the degree in $\overline{G}$ of each vertex is in $\{0,1,\ldots, (2n/3)^{1/3}-1\}$,  there are at least
$(\frac{3}{2})^{1/3}n^{2/3}$ vertices with the same degree in $\overline{G}$. Let $v$ be an arbitrary one of these vertices. Observe 
\[\Tri_{\overline{G}}(v) + e_{\overline{G}}\left(N_G(v),\overline{N_G(v)}\right)\le \binom{(2n/3)^{1/3}}{2}+(2n/3)^{2/3}<\left(\frac{3}{2}\right)^{1/3}n^{2/3}.\]
By the Pigeonhole Principle, we can find two distinct vertices which satisfy the hypotheses in  Lemma~\ref{comp}, so $G$ is not triangle-distinct. 

For (c), the lower bound follows directly from Theorem~\ref{thm-min-max-degree}(a). We now show the upper bound. 
Suppose for contradiction that $d > n - \sqrt{2n/3}$. 
Every vertex in $\overline{G}$ has degree $n - 1 - d < \sqrt{2n/3}$. 
Thus, for each vertex $v$, we have
\[\Tri_{\overline{G}}(v) + e_{\overline{G}}\left(N_G(v), \overline{N_G(v)}\right) < \binom{\sqrt{2n/3}}{2}+\left(\sqrt{2n/3}\right)^2 < n.\] 
By the Pigeonhole Principle, there are two distinct vertices which satisfy the conditions in Lemma~\ref{comp}, so $G$ is not triangle-distinct, a contradiction. 
\qed

%
%
%
Now we turn our attention to the number of edges in a triangle-distinct graph, but first we need a couple preliminary results.

\begin{lemma}\label{lemma:neighbor}
    Let $G$ be triangle-distinct with $n$ vertices. Let $c>0$ be a real number. 	If $e\left(\overline{G}\right) \leq cn$, then for any $1\leq k \leq n$ and any $0\leq t \leq k-1$, the number of vertices of degree $k-1$ in $\overline{G}$ which have at least $k-1-t$ common neighbors is at most $\sum_{i=0}^t (4cn)^{1-1/2^{i}}$.
\end{lemma}

\begin{proof}
Fix $1\leq k \leq n$ and let $V_k$ be the set of vertices of degree $k-1$ in $\overline{G}$. For each $0\leq t \leq k-1$, let $r_t$ be the maximum number of vertices $\{y_1, \ldots, y_{r_t}\}$ in $V_k$ for which $\left|\bigcap_{i=1}^{r_t} N_{\overline{G}}(y_i)\right| \geq (k-1)-t$.

We will use induction to prove $r_t\leq \sum_{i=0}^t (4cn)^{1-1/2^{i}}$. For the base case when $t=0$, notice that two vertices cannot share the same neighborhood because $G$ is assumed to be triangle-distinct. Therefore $r_0 \leq 1 = (4cn)^{0}$. 

For the inductive step, consider $1\leq t \leq k-1$. Suppose $\{y_1, y_2, \ldots, y_t\}$ are vertices in $V_k$ with $\left|\bigcap_{i=1}^{r_t} N_{\overline{G}}(y_i)\right| \geq (k-1)-t$. Let $U \subseteq \bigcap_{i=1}^{r_t} N_{\overline{G}}(y_i)$ with $|U| = (k-1)-t$. For each $1\leq i \leq r_t$, let $U_i = N_{\overline{G}}(y_i) \setminus U$. 

Because $G$ is triangle-distinct and for any $1\leq i<j \leq r_t$, $d_{\overline{G}}(y_i) = k-1 = d_{\overline{G}}(y_j)$, so Lemma~\ref{comp} implies 
\[t_{\overline{G}}(y_i) + e_{\overline{G}}\left(N_{G}(y_i), N_{\overline{G}}(y_i)\right)
\neq
t_{\overline{G}}(y_j) + e_{\overline{G}}\left(N_{G}(y_j), N_{\overline{G}}(y_j)\right).\]

For each $1\leq i\leq r_t$, define $p_i = e_{\overline{G}}(U, U_i)$ and $q_i = e\left({\overline{G}}[U_i]\right)$. Observe
\[t_{\overline{G}}(y_i) + e_{\overline{G}}\left(N_{G}(y_i), N_{\overline{G}}(y_i)\right) = 
\left(\sum_{u\in U_i} d_{\overline{G}}(u)\right) + \left(\sum_{u\in U} d_{\overline{G}}(u)\right)
-q_i -p_i - e\left(\overline{G}[U]\right) - (k-1).\]
Therefore, for each $1\leq i<j \leq r_t$,  
\[\left(\sum_{u\in U_i} d_{\overline{G}}(u)\right)-q_i -p_i 
\neq 
\left(\sum_{u\in U_j} d_{\overline{G}}(u)\right)-q_j -p_j.\]
Notice $\left(\sum_{u\in U_i} d_{\overline{G}}(u)\right)-q_i -p_i = e_{\overline{G}}\left( U_i, N_{G}(y_i) \right) + q_i + t \geq 0$.
Because these values are all distinct and non-negative, 
\[\sum_{i=1}^{r_t}\left( \left(\sum_{u\in U_i} d_{\overline{G}}(u)\right)-q_i -p_i \right)
\geq 
\sum_{i=1}^{r_t} (i-1) = \frac{1}{2}r_t(r_t-1)\geq \frac{1}{2}(r_t-1)^2.\]

Next, consider any $v\not\in U$ and $S_v=\{i\in \{1,\ldots, r_t\}: v\in U_i\}$. Observe that $U\cup\{v\}$ is a set of $(k-1)-(t-1)$ common neighbors for $\{y_i: i\in S_v\}$. Therefore, by the definition of $r_{t-1}$, we can conclude $|S_v| \leq r_{t-1}$. As a result, in the sum $\sum_{i=1}^{r_t} \left(\sum_{u\in U_i} d_{\overline{G}}(u)\right)$, the degree of each vertex is counted at most $r_{t-1}$ times. Therefore 
\[\sum_{i=1}^{r_t}\left( \left(\sum_{u\in U_i} d_{\overline{G}}(u)\right)-q_i-p_i \right)
\leq  \sum_{i=1}^{r_t}\sum_{u\in U_i} d_{\overline{G}}(u)
\leq 2r_{t-1}e\left(\overline{G}\right).\]
Combining the last two inequalities, together with the induction hypothesis and the assumption that $e(\overline{G}) \leq cn$, we find 
\begin{eqnarray*}
\frac{1}{2}(r_t-1)^2 & \leq & 2r_{t-1}e\left(\overline{G}\right) \leq 2cnr_{t-1}\leq 2cn\sum_{i=0}^{t-1} (4cn)^{1-1/2^{i}} = \frac{1}{2} \sum_{i=0}^{t-1} (4cn)^{2-1/2^{i}}\\
r_t-1 &\leq& \sqrt{\sum_{i=0}^{t-1} (4cn)^{2-1/2^{i}}} \leq \sum_{i=0}^{t-1} (4cn)^{1-1/2^{i+1}} = \sum_{i=1}^{t} (4cn)^{1-1/2^{i}} \\
r_t & \leq & \sum_{i=0}^{t} (4cn)^{1-1/2^{i}}.
\end{eqnarray*}

\end{proof}

\begin{corollary} \label{cor:deg}
	Let $G$ be triangle-distinct with $n$ vertices. Let $c>0$ be a real number.
	If $e(G) \geq {n \choose 2} - cn$, then for any $k \in \{1,\ldots, n\}$, the number of vertices of degree $n-k$ in $G$ is at most $k(4cn)^{1-1/2^{k-1}}$.
\end{corollary}

\begin{proof}
    Fix $k\in  \{1,\ldots, n\}$ and let $V_k$ be the set of vertices of degree $n-k$ in $G$. Note that	$V_k$ is also the set of vertices of degree $k - 1$ in $\overline{G}$. 

For each $0\leq t \leq k-1$, let $r_t$ be the maximum number of vertices $\{y_1, \ldots, y_{r_t}\}$ in $V_k$ for which $\left|\bigcap_{i=1}^{r_t} N_{\overline{G}}(y_i)\right| \geq (k-1)-t$. Observe $r_{k-1}=|V_{k}|$.

By Lemma~\ref{lemma:neighbor},
\[|V_k| = r_{k-1} \leq \sum_{i=0}^{k-1} (4cn)^{1-1/2^{i}} \leq \sum_{i=0}^{k-1} (4cn)^{1-1/2^{k-1}} = k(4cn)^{1-1/2^{k-1}}.\]
    
\end{proof}

Now let's use these results to bound the number of edges in a triangle-distinct graph.

\begin{theorem} \label{theorem:EdgeLow}
   Let $G$ be a triangle-distinct graph of order $n$. 
   Then 
   \[
     \frac{\sqrt{2}}{3}n^{3/2}\left(1-o(1)\right) <  e(G) \leq \binom{n}{2} - \omega(n). 
   \]
 \end{theorem}

\begin{proof}
   Let $G$ be a triangle-distinct graph with $n$ vertices.
   We'll first consider the lower bound.
   Notice that a vertex of degree $d$ is contained in at most $\binom{d}{2}$ triangles.
   Since $G$ is triangle-distinct, $G$ has at most $\binom{d}{2}+1$ vertices of degree at most $d$.
   Let $D$ be the largest integer such that $\binom{D}{2}+1 \leq n$. Then $\binom{D+1}{2} \geq n$, so $D > \sqrt{2n}-1$.
   
   Now let $v_1, v_2, \ldots, v_n$ be the vertices of $G$ in ascending degree order.
   For any non-negative $d$, there are at most $\binom{d}{2}+1$ vertices of degree at most $d$, so $\deg(v_i) > d$ for $i > \binom{d}{2}+1$.
   This gives us the lower bound
   \[\deg(v_i) \geq d \text{ for } \binom{d-1}{2} + 2 \leq i \leq \binom{d}{2} + 1.\]

   Therefore, the sum of degrees of the vertices in $G$ can be bounded as follows:
   \begin{eqnarray*}
   2e(G) 
   &=& \sum_{i=1}^{n} \deg(v_i)\\
   &\geq& \sum_{i=2}^{\binom{D}{2} + 1} \deg(v_i)\\
   &=& \sum_{d=2}^D  \sum_{i=\binom{d-1}{2} + 2}^{\binom{d}{2} + 1} \deg(v_i)\\
   &\geq& \sum_{d=2}^D  \sum_{i=\binom{d-1}{2} + 2}^{\binom{d}{2} + 1} d\\
   &=& \sum_{d=2}^D d\left( \binom{d}{2}  - \binom{d-1}{2}\right)\\
   &=& \sum_{d=2}^D d(d-1)\\
   &=& \frac{1}{3}(D+1)D(D-1)\\
   &>& \frac{1}{3}\left(\sqrt{2n}\right)\left(\sqrt{2n}-1\right)\left(\sqrt{2n}-2\right)\\
   &>& \frac{1}{3}\left(\sqrt{2}\left(n^{1/2}\right)-2\right)^3
   \end{eqnarray*}
This proves the lower bound. 

Next we turn our attention to the upper bound.
For each $1\leq i \leq n$, let $t_i$ be the number of vertices of degree $n-i$ in $G$. 
Suppose that $e(G) \geq {n \choose 2} - \frac{M}{2}n$ for some fixed $M>0$. Observe 
\begin{eqnarray*}
	e(G) & = & \frac{1}{2}\sum_{i = 1}^{n} t_i(n-i) \\
	&=& \frac{1}{2}\sum_{i = 1}^{M} t_i(n-i) + \frac{1}{2}\sum_{i = M+1}^{n} t_i(n-i)\\
	&\leq& \frac{1}{2}\sum_{i = 1}^{M} t_i(n-i) + \frac{1}{2}\left(\sum_{i=M+1}^n t_i\right)(n-M-1)\\
	&=& \frac{1}{2}\sum_{i = 1}^{M} t_i(n-i) + \frac{1}{2}\left(n-\sum_{i=1}^M t_i\right)(n-M-1)\\
	&=& \frac{1}{2}\sum_{i = 1}^{M} t_i(M-i+1) + \frac{1}{2}(n(n-1)-nM)\\
	&\leq & \frac{1}{2}\sum_{i = 1}^{M} i(2Mn)^{1-1/2^{i-1}}(M-i+1) + \frac{1}{2}(n(n-1)-nM)\quad \text{(by Corollary~\ref{cor:deg})}\\
	&\leq& \frac{1}{2}\sum_{i = 1}^{M} M(2Mn)^{1-1/2^{M-1}}(M) + \binom{n}{2} - \frac{M}{2}n \quad \text{(since $1 \leq i\leq M$)}\\
	&=& \frac{1}{2} M^3(2M)^{1-1/2^{M-1}} n^{1-1/2^{M-1}} + \binom{n}{2} - \frac{M}{2}n\\
	\end{eqnarray*}
 
So for any $M>0$, 
\[e(G) \geq {n \choose 2} - \frac{M}{2}n \quad\Longrightarrow\quad e(G) \leq {n \choose 2} - \frac{M}{2}n + f_{M}(n),\]
where 
\[f_{M}(n) = \left(\frac{1}{2} M^3(2M)^{1-1/2^{M-1}}\right) n^{1-1/2^{M-1}} > 0.\]
Thus it's impossible for $e(G) > {n \choose 2} - \frac{M}{2}n + f_{M}(n)$. And for fixed $M$, we asymptotically have $f_{M}(n) = o(n)$, so $e(G) \leq {n \choose 2} - \left(\frac{M}{2} - o(1)\right)n$. Since this it true for all $M$ no matter how large, $e(G) \leq \binom{n}{2} - \omega(n)$.

\end{proof}

 A planar graph on $n$ vertices has at most $3n-6$ edges. As a consequence of Theorem~\ref{theorem:EdgeLow}, 
 no triangle-distinct graph is planar; nor is its complement planar for sufficiently large $n$. 

\section{Open problems}

There is still much to learn about triangle-distinct graphs. For example, in Theorem~\ref{thm-min-max-degree} we show that if $G$ is $n$-vertex, $d$-regular, and triangle-distinct then $\sqrt{2n} < d \le n - \sqrt{2n/3}$. However, we have not found a regular triangle-distinct graph.

\begin{pro}
	Does there exist a regular graph that is triangle-distinct?
\end{pro}

 In Theorem~\ref{theorem:EdgeLow} we give bounds on the number of edges in a triangle-distinct graph of order $n$. However, we do not have examples showing these bounds are tight. 

\begin{pro}
Are the bounds on $e(G)$ for an $n$-vertex triangle-distinct graphs given in Theorem~\ref{theorem:EdgeLow} tight, and if not, how can we sharpen them? 
\end{pro}

\section{Acknowledgements}
This research was started at the Graduate Research Workshop in
Combinatorics (GRWC) held at Iowa State University in June 2015 and which was supported
by the grant NSF-DMS \#1500662.
The authors would also like to thank Sylvia Hobart, Bill Kay, and Brent McKain for helpful conversations and contributions in the early stages of this paper.


\begin{thebibliography}{100}


\bibitem{1} B. R. Nair and A. Vijayakumar, About triangles in a graph and its complement, {\it Discrete Math.} 131 (1994) 205-210. MR 95f:05058.

\bibitem{2} B. R. Nair and A. Vijayakumar, Strongly edge triangle regular graphs and a conjecture of Kotzig, {\it Discrete Math.} 158 (1996) 201-209. MR 97d:05283

\bibitem{personal} W. T. Trotter, personal communication.

\end{thebibliography}
\end{document}